\documentclass[preprint,11pt]{elsarticle}
\usepackage{amssymb}
\usepackage{amsfonts}
  \usepackage{color}
 \usepackage{bbm}
 \usepackage{mathrsfs,latexsym}
 \usepackage[all]{xy}
 \usepackage{amsthm,amsmath,amsfonts,amssymb}

\newtheorem{lem}{Lemma}[section]

\newtheorem{cor}{Corollary}[section]

\newtheorem{Acknw}{Acknowledgement}

\theoremstyle{definition}

\def\beq{\begin{equation}}
\def\eeq{\end{equation}}
\def\b{\begin}
\def\e{\end}

\begin{document}
\numberwithin{equation}{section}
\renewcommand{\theequation}{\arabic{section}.\arabic{equation}}

\begin{frontmatter}

\title{Potential theory in several quaternionic variables}

\author{Dongrui Wan}

\address{College of Mathematics and Computational Science, Shenzhen University, Room 412, Science and Technology Building, Shenzhen, 518060, P. R. China, Email: wandongrui@szu.edu.cn}

\begin{abstract}In this paper, we establish
the quaternionic versions of the potential description of various "small" sets related to the quaternionic plurisubharmonic functions in $\mathbb{H}^n$. We use the quaternionic capacity introduced in \cite{wan4} to characterize the $(-\infty)$-sets of plurisubharmonic functions, as the sets of the vanishing capacity. The latter requirement is also equivalent to the negligibility of the set. We also prove the Josefson's theorem on the equivalence of the locally and globally quaternionic polar sets in $\mathbb{H}^n$, following the method of Bedford-Taylor.
\end{abstract}


\end{frontmatter}

\section{Introduction}
The pluripotential theory, which is a non-linear complex counterpart of classical potential theory, has occupied an important place in mathematics. Although relatively young, the pluripotential theory has attracted considerable interest among analysts. The central part of the pluripotential theory is occupied by maximal plurisubharmonic functions and the generalized complex Monge-Amp\`{e}re operator. Decisive progress in this field has been made by Bedford and Taylor \cite{bed1980b,bed1980a,bed1990,bed}, Demailly \cite{Demailly1987,Demailly1991,demailly}, Cegrell \cite{cegrell1988,cegrell1998,cegrell2004,cegrell2009}, to mention only a few. Cegrell's book \cite{cegrell1988} provides an excellent in-depth study of capacities in $\mathbb{C}^n$. See also \cite{klimek} for detailed discussion on various type of small sets in $\mathbb{C}^n$.

In $n$-dimensional quaternionic space $\mathbb{H}^n$, little is known about the quaternionic pluripolar sets and the zero sets of the quaternionic capacities at present. The purpose of this paper is to give the potential-theoretic description of various "small" sets related to the quaternionic Monge-Amp\`{e}re operator.

Let $\Omega$ be an open set in $\mathbb{H}^n$. The quaternionic Monge-Amp\`{e}re operator is defined as the Moore determinant of the quaternionic Hessian of $u$:
 \begin{equation*}det(u)=det\left[\frac{\partial^2u}{\partial q_j\partial \bar{q}_k}(q)\right].\end{equation*}Alesker proved in \cite{alesker1} a quaternionic version of Chern-Levine-Nirenberg estimate and extended the definition of quaternionic Monge-Amp\`{e}re operator to continuous quaternionic plurisubharmonic functions. Since it is inconvenient to use the Moore determinant, the studying of the quaternionic Monge-Amp\`{e}re operator is much more difficult than the complex Monge-Amp\`{e}re operator.

To define the quaternionic Monge-Amp\`{e}re operator on general quaternionic manifolds, Alesker introduced in \cite{alesker2} an operator in terms of the Baston operator $\triangle$, which is the first operator of the quaternionic complex on quaternionic manifolds. The $n$-th power of this operator is exactly the quaternionic Monge-Amp\`{e}re operator when the manifold is flat. On the flat space $\mathbb{H}^n$, the Baston operator $\triangle$ is the first operator of $0$-Cauchy-Fueter complex:\begin{equation}\label{1.1}0\rightarrow C^\infty(\Omega,\mathbb{C})\xrightarrow{\triangle}C^\infty(\Omega,\wedge^2\mathbb{C}^{2n})
\xrightarrow{D}C^\infty(\Omega,\wedge^3\mathbb{C}^{2n})\rightarrow\cdots.
\end{equation}Wang \cite{Wang} wrote down explicitly each operator of the $k$-Cauchy-Fueter complex in terms of real variables.

Motivated by this, Wang and the author introduced in \cite{wan3} two first-order differential operators $d_0$ and $d_1$ acting on the quaternionic version of differential forms. And the second operator $D$ in (\ref{1.1}) can be written as $D:=\left(
                                                                          \begin{array}{c}
                                                                            d_0 \\
                                                                            d_1\\
                                                                          \end{array}
                                                                        \right)$.
The behavior of $d_0,d_1$ and $\triangle=d_0d_1$ is very similar to $\partial ,\overline{\partial}$ and $\partial\overline{\partial}$ in several complex variables. The quaternionic Monge-Amp\`{e}re operator can be defined as $(\triangle u)^n$ and has a simple explicit expression, which is much more convenient than the definition by using  Moore determinant.

By introducing the quaternionic version of differential forms, Wang and the author defined in \cite{wan4} the notions of closed positive forms and closed positive currents in the quaternionic case and our definition of closedness matches positivity well. We proved that $\triangle u$ is a closed positive $2$-current for any plurisubharmonic function $u$, and showed that when
functions $u_1,\ldots , u_k$ are locally bounded, $\triangle u_1\wedge\ldots\wedge\triangle u_k$ is a well defined closed positive current
and is continuous on decreasing sequences. In particular, the quaternionic Monge-Amp\`{e}re measure $(\triangle u)^n$ is well defined for locally bounded PSH function $u$, and is continuous on decreasing sequences converging to $u$.

Based on these observation, Zhang and the author established in \cite{wan4} several useful quaternionic versions of  results in the complex pluripotential theory, which play key roles in this paper. We showed that quasicontinuity, one of the most important properties of complex plurisubharmonic functions, holds also for quaternionic plurisubharmonic functions in $\mathbb{H}^n$. And we proved an equivalent characterization of the maximal PSH functions and the comparison theorems, which are connected to the uniqueness of Dirichlet problem of quaternionic Monge-Amp\`{e}re equations \cite{alesker4,jingyong}.

In this paper, we are concerned with the quaternionic capacities, and with the description of exceptional sets related to the quaternionic Monge-Amp\`{e}re operator in $\mathbb{H}^n$.

Now we introduce two types of exceptional sets in $\mathbb{H}^n$. A set $E\subseteq\mathbb{H}^n$ is said to be \emph{locally quaternionic polar} (locally Q-polar, for short) if for each point $a\in E$, there is a neighborhood  $B(a,r)$ and a function $u\in PSH(B(a,r))$ such that $u|_{E\cap B(a,r)}=-\infty$. And a set $E\subseteq\Omega$ in $\mathbb{H}^n$ is said to be \emph{globally quaternionic polar} (globally Q-polar, for short) in $\Omega$ if there exists a function $u\in PSH(\Omega)$ such that $E\subseteq\{u=-\infty\}$.

We show the proof of Josefson's theorem on the equivalence of the locally and globally Q-polar sets in $\mathbb{H}^n$, following the proof in pluripotential theory on $\mathbb{C}^n$ given in \cite{bed}. The original proof was given by Josefson \cite{josefson} basing on complicated estimates for polynomials.
\b{thm}\label{josefson}If $P\subseteq\mathbb{H}^n$ is locally Q-polar, there exists $v\in PSH(\mathbb{H}^n)$ with $P\subseteq\{v=-\infty\}$, i.e. $P$ is globally Q-polar in $\mathbb{H}^n$.
\e{thm}

We consider also the so-called \emph{negligible sets}, which are those of the form
\beq \label{neg}N=\{q\in\Omega,~u(q)<u^*(q)\},\eeq
where $u=\sup_\alpha u_\alpha$ is the upper envelope of a family of functions $(u_\alpha)\subseteq PSH(\Omega)$ which are locally bounded from above in $\Omega$, and $u^*$ is the upper semicontinuous regularization of $u$, i.e. $u^*(q)=\limsup_{q'\rightarrow q}u(q')$, for $q\in \Omega.$ We show in Theorem \ref{t1} that the negligible sets are precisely the Q-polar sets.

In order to prove the quasicontinuity theorem (Lemma \ref{t2.4} in Section 2), the author introduced in \cite{wan4} the quaternionic capacities for quaternionic plurisubharmonic functions. These capacities are defined in the same way as the capacities introduced by Bedford and Taylor in \cite{bed} for plurisubharmonic functions in $\mathbb{C}^n$.

Let $\Omega$ be an open set of $\mathbb{H}^n$, and let $K\subseteq\Omega$ be a compact set. The \emph{(relative) quaternionic capacity} of $K$ in $\Omega$ is defined by
\begin{equation}\label{capacity}C(K,\Omega)=\sup\left\{\int_{K}(\triangle
u)^n:u\in
PSH(\Omega),0\leq u\leq1\right\}.
\end{equation}
For any set $E\subseteq\Omega$, the\emph{ outer capacity} of $E$ is defined by
\beq\label{cap*}C^*(E,\Omega)=\inf\{C(\omega,\Omega):\omega~is~open,~~ E\subseteq\omega\subseteq\Omega\}.
\eeq
The main result of this paper is the following theorem giving the characterization of exceptional sets in terms of the outer capacity.

\b{thm}\label{t1}Let $\Omega\subseteq\mathbb{H}^n$ and $E\subseteq\Omega$. The following three statements are equivalent:\\
(1). $E$ is Q-polar;\\
(2). $E$ is negligible;\\
(3). $C^*(E,\Omega)=0$.\\
In particular, if $\Omega$ is strongly pseudoconvex smooth open set in $\mathbb{H}^n$ and $E\Subset\Omega$, then each of (1)-(3) is equivalent to (4): $u_E^*=0$.
\e{thm}
Here $u_E^*$ is the upper semicontinuous  regularization of the \emph{relative extremal function} $ u_E$ defined as ($E\Subset\Omega$)
\beq\label{u*} u_E(q)= u_{E,\Omega}(q)=\sup\{u(q):u\in PSH(\Omega),u\leq0,u|_E\leq -1\}, q\in \Omega.\eeq
The function $u_E^*$ is a powerful tool to give the connection between the outer capacities and the Q-polar sets.
\par
Finally we prove that the outer capacity $C^*(\cdot,\Omega)$ is a generalized capacity in the sense of Choquet (Theorem \ref{t3.1} in Section 3).

Although we strongly use ideas of Bedford and Taylor \cite{bed} from the pluripotential theory in $\mathbb{C}^n$, our potential results for the quaternionic Monge-Amp\`{e}re operator are completely new. The theory of quaternionic closed positive currents established recently in \cite{wan3} allows us to treat the quaternionic Monge-Amp\`{e}re operator as an operator of divergence form, and so we can integrate by parts. Since this can avoid the inconvenience in using Moore determinant, we established several useful quaternionic versions of results in the complex pluripotential theory in \cite{wan4}. All these preparation play key roles in this paper.

\section{Preliminaries on quaternionic Monge-Amp\`{e}re measure}\label{sec2}
Recall that an upper semicontinuous function $u$ on
$\mathbb{H}^{n}$ is said to be \emph{quaternionic plurisubharmonic} if $u$ is
subharmonic on each right quaternionic line. Denote by $PSH$ the class
of all quaternionic plurisubharmonic functions (cf.
\cite{alesker1,alesker4,alesker3,alesker2} for more information about quaternionic plurisubharmonic functions).

As in \cite{wan3}, we use the conjugate embedding
\begin{equation*}\label{2.2}\begin{aligned}\tau:\mathbb{H}^{n}\cong\mathbb{R}^{4n}&\hookrightarrow\mathbb{C}^{2n\times2},\\ (q_0,\ldots,q_{n-1})&\mapsto \textbf{z}=(z^{j\alpha})\in\mathbb{C}^{2n\times2},
\end{aligned}\end{equation*}
$q_j=x_{4j}+\textbf{i}x_{4j+1}+\textbf{j}x_{4j+2}+\textbf{k}x_{4j+3}$, $j=0,1,\ldots,2n-1, ~\alpha=0 ,1 ,$ with
\begin{equation}\label{2.31}\left(
                             \begin{array}{cc}
                               z^{00 } & z^{01 } \\
                               z^{10 } & z^{11 } \\
                               \vdots&\vdots\\
                               z^{(2l)0 } & z^{(2l)1 } \\
                               z^{(2l+1)0 } & z^{(2l+1)1 } \\
                               \vdots&\vdots\\
                               z^{(2n-2)0 } & z^{(2n-2)1 } \\
                               z^{(2n-1)0 } & z^{(2n-1)1 } \\
                             \end{array}
                           \right):=\left(
                                      \begin{array}{cc}
                                      x_{0}-\textbf{i}x_{1} & -x_{2}+\textbf{i}x_{3} \\
                                        x_{2}+\textbf{i}x_{3} & x_{0}+\textbf{i}x_{1} \\
                                        \vdots&\vdots\\
                                         x_{4l}-\textbf{i}x_{4l+1} & -x_{4l+2}+\textbf{i}x_{4l+3} \\
                                        x_{4l+2}+\textbf{i}x_{4l+3} & x_{4l}+\textbf{i}x_{4l+1} \\
                                        \vdots&\vdots\\
                                        x_{4n-4}-\textbf{i}x_{4n-3} & -x_{4n-2}+\textbf{i}x_{4n-1} \\
                                        x_{4n-2}+\textbf{i}x_{4n-1} & x_{4n-4}+\textbf{i}x_{4n-3} \\
                                      \end{array}
                                    \right).
\end{equation} Pulling back to the quaternionic
space $\mathbb{H}^n\cong\mathbb{R}^{4n}$ by the embedding
(\ref{2.31}), we define on $\mathbb{R}^{4n}$ first-order differential operators $\nabla_{j\alpha}$ as following:
\begin{equation}\label{2.4}\left(
                             \begin{array}{cc}
                               \nabla_{00 } & \nabla_{01 } \\
                               \nabla_{10 } & \nabla_{11 } \\
                               \vdots&\vdots\\
                               \nabla_{(2l)0 } & \nabla_{(2l)1 } \\
                               \nabla_{(2l+1)0 } & \nabla_{(2l+1)1 } \\
                               \vdots&\vdots\\
                               \nabla_{(2n-2)0 } & \nabla_{(2n-2)1 } \\
                               \nabla_{(2n-1)0 } & \nabla_{(2n-1)1 } \\
                             \end{array}
                           \right):=\left(
                                      \begin{array}{cc}
                                      \partial_{x_{0}}+\textbf{i}\partial_{x_{1}} & -\partial_{x_{2}}-\textbf{i}\partial_{x_{3}} \\
                                        \partial_{x_{2}}-\textbf{i}\partial_{x_{3}} & \partial_{x_{0}}-\textbf{i}\partial_{x_{1}} \\
                                         \vdots&\vdots\\
                                         \partial_{x_{4l}}+\textbf{i}\partial_{x_{4l+1}} & -\partial_{x_{4l+2}}-\textbf{i}\partial_{x_{4l+3}} \\
                                        \partial_{x_{4l+2}}-\textbf{i}\partial_{x_{4l+3}} & \partial_{x_{4l}}-\textbf{i}\partial_{x_{4l+1}} \\
                                        \vdots&\vdots\\
                                        \partial_{x_{4n-4}}+\textbf{i}\partial_{x_{4n-3}} & -\partial_{x_{4n-2}}-\textbf{i}\partial_{x_{4n-1}} \\
                                        \partial_{x_{4n-2}}-\textbf{i}\partial_{x_{4n-1}} & \partial_{x_{4n-4}}-\textbf{i}\partial_{x_{4n-3}} \\
                                      \end{array}
                                    \right).
\end{equation}

Let $\wedge^{2k}\mathbb{C}^{2n}$ be the complex exterior algebra generated by $\mathbb{C}^{2n}$, $0\leq k\leq n$. Fix a basis
$\{\omega^0,\omega^1,\ldots$, $\omega^{2n-1}\}$ of $\mathbb{C}^{2n}$. Let $\Omega$ be a domain in $\mathbb{R}^{4n}$. We
define $d_0,d_1:C_0^\infty(\Omega,\wedge^{p}\mathbb{C}^{2n})\rightarrow C_0^\infty(\Omega,\wedge^{p+1}\mathbb{C}^{2n})$ by \begin{equation*}\begin{aligned}\label{2.228}&d_0F=\sum_{k,I}\nabla_{k0 }f_{I}~\omega^k\wedge\omega^I,\\
&d_1F=\sum_{k,I}\nabla_{k1 }f_{I}~\omega^k\wedge\omega^I,\\
&\triangle
F=d_0d_1F,
\end{aligned}\end{equation*}for $F=\sum_{I}f_{I}\omega^I\in C_0^\infty(\Omega,\wedge^{p}\mathbb{C}^{2n})$,  where the multi-index
$I=(i_1,\ldots,i_{p})$ and
$\omega^I:=\omega^{i_1}\wedge\ldots\wedge\omega^{i_{p}}$. Although $d_0,d_1$ are not exterior differential, their behavior is similar to the exterior differential: $d_0d_1=-d_1d_0$; $d_0^2=d_1^2=0$; for $F\in C_0^\infty(\Omega,\wedge^{p}\mathbb{C}^{2n})$, $G\in C_0^\infty(\Omega,\wedge^{q}\mathbb{C}^{2n})$, we have\begin{equation}\label{da}d_\alpha(F\wedge G)=d_\alpha F\wedge G+(-1)^{p}F\wedge d_\alpha G, \quad\alpha=0,1,\qquad d_0\triangle=d_1\triangle=0.\end{equation}(\ref{1.1}) is a complex since $D\triangle=0$.

For $u_1,\ldots,
u_n\in C^2$, it follows easily from (\ref{da}) that $\triangle u_1\wedge\ldots\wedge\triangle u_n$ satisfies the following remarkable identities:\begin{equation*}\begin{aligned}&\triangle u_1\wedge \triangle
u_2\wedge\ldots\wedge\triangle u_n=d_0(d_1u_1\wedge \triangle
u_2\wedge\ldots\wedge\triangle u_n)\\&=-d_1(d_0u_1\wedge \triangle
u_2\wedge\ldots\wedge\triangle u_n)=d_0d_1(u_1\triangle
u_2\wedge\ldots\wedge\triangle u_n)\\&=\triangle (u_1
\triangle u_2\wedge\ldots\wedge\triangle u_n).
\end{aligned}\end{equation*}

We
define for $u\in C^2$, $\triangle_{ij}u:=\frac{1}{2}(\nabla_{i0 }\nabla_{j1 }u-\nabla_{i1 }\nabla_{j0 }u).$
Then for $u_1,\ldots,u_n\in
C^2$,
\begin{equation*}\label{2.11}\begin{aligned}\triangle
u_1\wedge\ldots\wedge\triangle
u_n&=\sum_{i_1,j_1,\ldots}\triangle_{i_1j_1}u_1\ldots\triangle_{i_nj_n}u_n~\omega^{i_1}\wedge
\omega^{j_1}\wedge\ldots\wedge \omega^{i_n}\wedge
\omega^{j_n}\\&=\sum_{i_1,j_1,\ldots}\delta^{i_1j_1\ldots
i_nj_n}_{01\ldots(2n-1)}\triangle_{i_1j_1}u_1\ldots\triangle_{i_nj_n}u_n~\Omega_{2n},
\end{aligned}\end{equation*}where $\Omega_{2n}$ is defined as\begin{equation}\label{2.21}\Omega_{2n}:=\omega^0\wedge
\omega^1\wedge\ldots\wedge\omega^{2n-2}\wedge
\omega^{2n-1},\end{equation}and $\delta^{i_1j_1\ldots
i_nj_n}_{01\ldots(2n-1)}:=$ the sign of the permutation from $(i_1,j_1,\ldots
i_n,j_n)$ to\\ $(0,1,\ldots,2n-1),$ if $\{i_1,j_1,\ldots,
i_n,j_n\}=\{0,1,\ldots,2n-1\}$; otherwise,\\ $\delta^{i_1j_1\ldots
i_nj_n}_{01\ldots(2n-1)}=0$. In particular, when
$u_1=\ldots=u_n=u$, $\triangle u_1\wedge\ldots\wedge\triangle u_n$
coincides with $(\triangle u)^n:=\wedge^n\triangle u$.

Although a $2n$-form is not an authentic differential form and we cannot integrate it, we can define $\int_\Omega F:=\int_\Omega f dV,$ if we write $F=f~\Omega_{2n}\in L^1(\Omega,\wedge^{2n}\mathbb{C}^{2n})$,
where $dV$ is the Lesbesgue measure and $\Omega_{2n}$ is given by (\ref{2.21}). In particular, if $F$ is positive $2n$-form, then $\int_\Omega F
\geq0$. For a $2n$-current $F=\mu~\Omega_{2n}$ with coefficient to be measure $\mu$, define
\begin{equation*}\int_\Omega F:=\int_\Omega \mu.\end{equation*}

We proved that $\triangle u$ is a closed positive $2$-current for any $u\in PSH(\Omega)$.
Inductively, for $u_1,\ldots,u_p\in PSH\cap
L_{loc}^\infty(\Omega)$, we showed that
\begin{equation}\label{3.111}\triangle
u_1\wedge\ldots\wedge\triangle u_p :=\triangle(u_1\triangle
u_2\ldots\wedge\triangle u_p )\end{equation} is a closed
positive $2p$-current. In particular, for $u_1,\ldots,u_n\in PSH\cap
L_{loc}^{\infty} (\Omega)$, $\triangle u_1\wedge\ldots\wedge\triangle u_n=\mu\Omega_{2n}$ for a well defined positive Radon
measure $\mu$. See \cite{wan3} for the detailed information about the closed positive currents in $\mathbb{H}^n$.\par

\begin{lem}\label{l2.3}(Theorem 3.1 in \cite{wan3}) Let
$v^1,\ldots,v^k\in PSH\cap L_{loc}^{\infty}(\Omega)$ and let
$\{v_j^1\}_{j\in\mathbb{N}},$ $\ldots,\{v_j^k\}_{j\in\mathbb{N}}$ be
decreasing sequences of $PSH$ functions in $\Omega$ such that $\lim_{j\rightarrow\infty}v_j^t=v^t$ pointwisely in $\Omega$
for each $t $. Then the currents $\triangle
v_j^1\wedge\ldots\wedge\triangle v_j^k $ converge weakly to
$\triangle v^1\wedge\ldots\wedge\triangle v^k $ as
$j\rightarrow\infty$.
\end{lem}
Alesker gave a quaternionic version of Chern-Levine-Nirenberg estimate in Proposition 6.3 in \cite{alesker2}.  And we gave  an elementary and simpler proof in \cite{wan3}.

\begin{lem}\label{p3.88}$($Chern-Levine-Nirenberg type estimate, see Proposition 3.10 in \cite{wan3} $)$ Let
$\Omega$ be a domain in $\mathbb{H}^n$. Let $K,L$ be compact subsets
of $\Omega$ such that $L$ is contained in the interior of $K$. Then
there exists a constant $C$ depending only on $K,L$ such that for
any $u_1,\ldots u_n\in PSH\cap C^2(\Omega)$, one
has
\begin{equation}\label{3.14}\int_L\triangle
u_1\wedge\ldots\wedge\triangle u_n \leq
C\prod_{i=1}^n\|u_i\|_{L^{\infty}(K)}.
\end{equation} \end{lem}This estimate also holds for any $u_1,\ldots u_n\in PSH\cap L_{loc}^\infty(\Omega)$.

The analogous classical results for subharmonic functions also hold for the quaternionic plurisubharmonic functions. We list these properties here without proofs; all of them can be derived from the subharmonic case (see Chapter 2 in \cite{klimek}).
\b{pro}\label{p2.1} Let $\Omega$ be an open subset of $\mathbb{H}^n$.\\
(1). The family $PSH(\Omega)$ is a convex cone, i.e. if $\alpha,\beta$ are non-negative numbers and $u,v\in PSH(\Omega)$, then $\alpha u+\beta v\in PSH(\Omega)$; and $\max\{u,v\}\in PSH(\Omega)$.\\
(2). If $\Omega$ is connected and $\{u_j\}\subseteq PSH(\Omega)$ is a decreasing sequence, then $u=\lim_{j\rightarrow\infty}u_j\in PSH(\Omega)$ or $u\equiv-\infty$.\\
(3). Let $\{u_\alpha\}_{\alpha\in A}\subseteq PSH(\Omega)$ be such that its upper envelope $u=\sup_{\alpha\in A}u_\alpha$ is locally bounded above. Then the upper semicontinuous regularization $u^*\in PSH(\Omega)$.\\
(4). Let $\omega$ be a non-empty proper open subset of $\Omega$, $u\in PSH(\Omega),v\in PSH(\omega)$, and $\limsup_{q\rightarrow\zeta}v(q)\leq u(\zeta)$ for each $\zeta\in \partial\omega\cap\Omega$, then
$$w:=\left\{
      \begin{array}{ll}
        \max\{u,v\}, & \text{in}~\omega \\
        u, &  \text{in}~\Omega\backslash\omega
      \end{array}
    \right. \quad\in PSH(\Omega).$$\\
(5). Let $F$ be a closed subset of $\Omega$ of the form $F=\{q\in \Omega,v(q)=-\infty\},$ where $v\in PSH(\Omega)$. If $u\in PSH(\Omega\backslash F)$ is bounded above, then  $$\widetilde{u}(q):=\left\{
                                                                                     \begin{array}{ll}
                                                                                       u(q), & q\in \Omega\backslash F \\
                                                                                       \limsup_{q'\notin F,q'\rightarrow q }u(q'), & q\in F
                                                                                     \end{array}
                                                                                   \right.  \quad\in PSH(\Omega).$$
\e{pro}

\begin{lem}\label{t2.4}(quasicontinuity theorem, see Theorem 1.1 in \cite{wan4}) Let $\Omega$ be an open subset of $\mathbb{H}^n$ and let $u$ be a locally bounded PSH function. Then for each $\varepsilon>0$, there exists an open subset $\omega$ of $\Omega$ such that $C(\omega)<\varepsilon$ and $u$ is continuous on $\Omega\backslash\omega$.
\end{lem}

After appropriate results are proved in Section 3 later, it will be clear in Remark 31 that the assumption of locally boundedness of the function $u$ is superfluous.

We have already seen in Lemma \ref{l2.3} that the quaternionic Monge-Amp\`{e}re operator is continuous on decreasing sequences of locally bounded PSH functions. It turns out that it also behaves well on increasing sequences just as the complex Monge-Amp\`{e}re operator. See \cite{bed} for the analogue result in the complex case.

\begin{lem}\label{l2}(Proposition 4.1 in \cite{wan4}) Let $\{u_j\}_{j\in \mathbb{N}}$ be a sequence in $PSH\cap L_{loc}^\infty(\Omega)$ that increases to $u\in PSH\cap L_{loc}^\infty(\Omega)$ almost everywhere in $\Omega$ $($with respect to Lebesgue measure$)$. Then the currents $(\triangle u_j)^n$ converge weakly to $(\triangle u)^n$ as $j\rightarrow\infty$.
\end{lem}

By using the quasicontinuity theorem and the convergence result, we proved in \cite{wan4} the following comparison theorem, which will be a useful tool in this paper.
This comparison result implies the minimum principle results in \cite{alesker1}, which are essential to the uniqueness of Dirichlet problem of quaternionic Monge-Amp\`{e}re equations (cf. \cite{alesker4,jingyong}).

\begin{lem}\label{t2.2}(See Theorem 1.2 in \cite{wan4}) Let $u,v\in PSH\cap L_{loc}^\infty(\Omega)$. If for any $\zeta\in\partial\Omega$,
$$\liminf_{\zeta\leftarrow q\in\Omega}~~(u(q)-v(q))\geq0,$$
then
\begin{equation}\label{5.1}\int_{\{u<v\}}(\triangle v)^n\leq\int_{\{u<v\}}(\triangle u)^n.
\end{equation}
\end{lem}
\section{Capacity and description of exceptional sets}

Let $\Omega$ be an open set of $\mathbb{H}^n$, and let $K\subseteq\Omega$ be a compact set. The relative quaternionic capacity $C(K,\Omega)$ is given by (\ref{capacity}).
For an arbitrary set $E\subseteq\Omega$, the standard construction of \emph{inner capacity} is
\beq\label{capacity2}C(E,\Omega)=\sup\{C(K,\Omega):K~is~compact,~~K\subseteq E\}.
\eeq
The Chern-Levine-Nirenberg type estimate (Lemma \ref{p3.88}) shows that $C(E,\Omega)<+\infty$ for any $E\subseteq\Omega$. We give some elementary properties of capacity.

\b{pro}\label{p1.1}(1). If $E\subseteq\Omega$ is Borel, then
$$ C(E,\Omega)=\sup\left\{\int_{E}(\triangle
u)^n:u\in
PSH(\Omega),0\leq u\leq1\right\}.
$$
(2). If $E_1\subseteq E_2\subseteq\Omega$, then $C(E_1,\Omega)\leq C(E_2,\Omega)$.\\
(3). If $E\subseteq \Omega_1\subseteq\Omega_2$, then $C(E,\Omega_1)\geq C(E,\Omega_2)$.\\
(4). If $E_1,E_2,\ldots\subseteq\Omega$, then $$C(\bigcup_{j=1}^\infty E_j,\Omega)\leq \sum_{j=1}^\infty C(E_j,\Omega).$$
(5). If $E_1\subseteq E_2\subseteq\ldots\subseteq\Omega$, then $$C(\bigcup_{j=1}^\infty E_j,\Omega)= \lim_{j\rightarrow\infty} C(E_j,\Omega).$$
(6). If $\Omega_1\subseteq\Omega_2\subseteq H^n$ and $\omega\Subset\Omega_1$, then there exists a constant $A>0$ such that for all Borel subsets $E\subseteq\omega$ we have $C(E,\Omega_1)\leq A~C(E,\Omega_2)$.
\e{pro}
\b{proof}   Note that (1)-(5) are the direct consequences of the definition (\ref{capacity}) and (\ref{capacity2}), we give the proof of property (6) here. Without loss of generality we may suppose that $\omega\Subset\Omega_1$ are concentric balls, say $\Omega_1=B(0,r)$ and $\omega=B(0,r-\varepsilon)$, $\varepsilon>0$. For each $u\in PSH(\Omega_1)$ with $0\leq u\leq1$, define
$$ \widetilde{u}(q)=\left\{
                      \begin{array}{ll}
                        \max\{u(q),\lambda(\|q\|^2-r^2)+2\}, & on~\Omega_1, \\
                        \lambda(\|q\|^2-r^2)+2, &on~\Omega_2\backslash\Omega_1.
                      \end{array}
                    \right.$$
Take a constant $\lambda$ sufficiently large such that $\lambda((r-\varepsilon)^2-r^2)\leq -2$. Then $\widetilde{u}=u$ on $\omega$, and $u(q)<\lambda(\|q\|^2-r^2)+2$ on $\partial\Omega_1$. By Proposition \ref{p2.1} we have $\widetilde{u}\in PSH(\Omega_2)$. Since $0\leq\widetilde{u}\leq M$ for some constant $M>0$, $0\leq\frac{\widetilde{u}}{M}\leq1$. For
any Borel subset $E\subseteq\omega$ we have
$$\int_E (\triangle u)^n=\int_E (\triangle \widetilde{u})^n\leq M^n C(E,\Omega_2).$$
Hence $C(E,\Omega_1)\leq M^n C(E,\Omega_2)$.
\e{proof}

\b{lem}\label{l14.4}Let $K$ be a compact subset of $\Omega$ and $\omega\Subset\Omega$ a neighborhood of $K$. There is a constant $A>0$ such that for each $v
\in PSH(\Omega)$, $$C(K\cap\{v<-m\},\Omega) \leq A\|v\|_{L^1(\overline{\omega})}\cdot\frac{1}{m}.$$
\e{lem}
\b{proof} For each $u\in PSH(\Omega)$, $0\leq u\leq 1$, by Chern-Levine-Nirenberg type estimate (Lemma \ref{p3.88}) we have $$\int_{K\cap\{v<-m\}}(\triangle u)^n\leq \frac{1}{m}\int_K|v|(\triangle u)^n \leq \frac{1}{m}C_{K,\overline{\omega}}\|v\|_{L^1(\overline{\omega})}.$$
\e{proof}

  Note that the outer capacity $C^*(\cdot,\Omega)$ given by (\ref{cap*}) also satisfies properties (2),(3),(4) in Proposition \ref{p1.1} above. $C(E,\Omega)\leq C^*(E,\Omega)$ for any $E\subseteq\Omega$. And it follows from the definition that $C(E,\Omega)=C^*(E,\Omega)$ for all open sets $E\subseteq\Omega$.

\b{cor}\label{c2.1} If $P$ is globally Q-polar in $\Omega$, then $C^*(P,\Omega)=0$.
\e{cor}
\proof By definition of globally Q-polar, assume that $P\subseteq\{v=-\infty\}$ with $v\in PSH(\Omega)$. Let $\Omega=\bigcup_{j\geq1}\Omega_j$ with $\Omega_j\Subset\Omega$. By Lemma \ref{l14.4} there is an open set $G_j=\Omega_j\cap\{v<-m_j\}$ with $C(G_j,\Omega)<\varepsilon2^{-j}$. So we have $\{v=-\infty\}\subseteq G=\bigcup G_j$ and $C(G,\Omega)<\varepsilon$.
\endproof

Let $\{u_\alpha\}$ be a family of PSH functions in $\Omega$ which is locally bounded from above. Then the function $u=\sup_\alpha u_\alpha$
is not in general $PSH$ or even upper semicontinuous. But its upper semicontinuous regularization $$u^*(q)=\limsup_{q'\rightarrow q}u(q')\geq u(q),~~~q\in \Omega$$ is PSH by Proposition \ref{p2.1}(3). A set of the form $$N=\{q\in \Omega: u(q)<u^*(q)\}$$ is called negligible. $u^*=u$ almost everywhere in $\Omega$ by the well known result for the subharmonic functions. So the Lebesgue measure of any negligible set $N$ is zero.

\b{pro}\label{p13.2} If $\Omega\subseteq\mathbb{H}^n$ is strongly pseudoconvex smooth open, then each globally Q-polar set $P=\{v=-\infty\}$ is negligible.
\e{pro}
\b{proof} Let $w$ be a smooth PSH function in $\Omega$, $w\geq v$. Denote $u_\alpha=(1-\alpha)v+\alpha w$, $\alpha\in (0,1)$. Then $u_\alpha$ is increasing in $\alpha$ and $$u:=\sup_\alpha u_\alpha=\left\{
                                \begin{array}{ll}
                                  -\infty, & \{v=-\infty\}, \\
                                  w, & \{v>-\infty\}
                                \end{array}
                              \right.$$ So $u^*=w$ in $\Omega$. $P=\{v=-\infty\}$ is exactly the set $\{u<u^*\}$.
\e{proof}

The main tool in the proof of our main result (Theorem \ref{t1})  is the relative extremal function $u_E$ defined by (\ref{u*}).
 Its upper semicontinuous regularization $u_E^*\in PSH(\Omega)$ (by Proposition \ref{p2.1}) and $-1\leq u_E^*\leq0$ in $\Omega$. If $\Omega$ is strongly pseudoconvex, then $u_E^*(q)\rightarrow0$ as $q\rightarrow\partial\Omega$. This function is defined in the same way as the relative extremal function given by Demailly \cite{Demailly1991} for the complex case, see also \cite{klimek} for detailed discussion. In the complex case, $u_E^*$ is sometimes called the PSH measure of $E$ relative to $\Omega$ \cite{sadu} or the regularized relative extremal function.

\b{lem}\label{l1} Fix a ball $\overline{B}\subseteq\Omega$. For any $g\in PSH\cap L_{loc}^\infty(\Omega)$, there exists a PSH function $\widetilde{g}$ such that $\widetilde{g}\geq g$ on $\Omega$ and $\widetilde{g}=g$ on $\Omega\backslash B$ and $(\triangle \widetilde{g})^n=0$ on $B$.
\e{lem}
\proof First consider the case when $g$ is continuous in $\Omega$. We shall use the Perron-Bremermann function defined by Alesker (Section 6 in \cite{alesker4}): $$u=\sup\{v:v~\text{is~finite~PSH~on~}B,\limsup_{q\rightarrow\zeta}v(q)\leq g(\zeta), ~~\forall~\zeta\in\partial B\}. $$
Then by Theorem 6.1 in \cite{alesker4}, $u\in PSH(B)$ is continuous on $\overline{B}$ and $u=g$ on $\partial B$. Note that $u$ is maximal in $B$, by Theorem 1.3 in \cite{wan4}, we have $(\triangle u)^n=0$ on $B$. Let $$\widetilde{g}=\left\{
                                                                                                                                                  \begin{array}{ll}
                                                                                                                                                    u, & \hbox{on}~~ B, \\
                                                                                                                                                    g, & \hbox{on}~~\Omega\backslash B,
                                                                                                                                                  \end{array}
                                                                                                                                                \right.$$ then $\widetilde{g}\geq g$ on $B$. Since $\widetilde{g}$ is the decreasing limit of PSH functions $$g_k=\left\{
                                                                                                                                                  \begin{array}{ll}
                                                                                                                                                    \max\{u,g+\frac{1}{k}\}, & \hbox{on}~~ B, \\
                                                                                                                                                   g+\frac{1}{k}, & \hbox{on}~~\Omega\backslash B,
                                                                                                                                                  \end{array}
                                                                                                                                                \right.$$ we get $\widetilde{g}\in PSH(\Omega)$ by Proposition \ref{p2.1}.

For an arbitrary function $g\in PSH\cap L_{loc}^\infty(\Omega)$, its regularization $g_{l}:=g*\rho_{\frac{1}{l}}\searrow g$. The function $\widetilde{g}:=\lim_{l\rightarrow+\infty}\widetilde{g}_{l}$ has all required properties.
\endproof
\b{lem}\label{choquet} (Choquet's lemma) Every family $(u_\alpha)$ has a countable subfamily $(u_{\alpha(j)})$ whose upper envelope $v$ satisfies $v\leq u\leq u^*=v^*$, where $u$ is the upper envelope of $(u_\alpha)$.
\e{lem}

\b{pro}\label{p13.5} Let $\Omega$ be an open set in $\mathbb{H}^n$ and let $K\subseteq \Omega$ be compact, then $(\triangle u_K^*)^n=0$ on $\Omega\backslash K$, that is, $u_K^*$ is maximal in $\Omega\backslash K$.
\e{pro}
\proof By Lemma \ref{choquet}, there exists a sequence $\{v_j\}\subseteq PSH(\Omega)$ such that $v_j\leq0$ on $\Omega$, $v_j\leq-1$ on $K$ and $v^*=u_K^*$. Replace $v_j$ by $\max\{-1,v_1,\ldots,v_j\}$ then we can assume that $\{v_j\}$ is increasing and $v_j\geq-1$ for all $j$.

Fix a ball $B\subseteq\Omega\backslash K$. Let $\widetilde{v_j}$ be as in Lemma \ref{l1}. We have $\widetilde{v_j}\leq0$ on $\Omega$,  $\widetilde{v_j}\leq-1$ on $K$. Then $v_j\leq \widetilde{v_j}\leq u_K$ and $\widetilde{v}=\lim_j\widetilde{v_j}$ such that $v^*=\widetilde{v}^*=u_K^*$ and $\lim \widetilde{v_j}=\lim v_j=u_K^*$ a.e. in $\Omega$. Since $(\triangle \widetilde{v_j})^n=0$ on $B$, by Lemma \ref{l2}, $(\triangle u_K^*)^n=0$ on $B$. And $B$ is taken arbitrary, it follows $(\triangle u_K^*)^n=0$ on $\Omega\backslash K$.
 \endproof

\b{cor}\label{c13.5} For arbitrary sets $E\Subset\Omega$, we have the following properties of the regularized relative extremal functions $u_E^*$:\\
(1). If $E_1\subseteq E_2\subseteq\Omega_1\subseteq\Omega_2$, then $u_{E_1,\Omega_1}^*\geq u_{E_2,\Omega_1}^*\geq u_{E_2,\Omega_2}^*$.\\
(2). $u_E^*=u_E=-1$ on $E^0$ and $(\triangle u_E^*)^n=0$ on $\Omega\backslash \overline{E}$; so $(\triangle u_E^*)^n$ is supported by $\partial E$.
\e{cor}
\b{proof} (1) is obvious. And it follows directly from the definition that $u_E^*=u_E=-1$ on $E^0$, hence $(\triangle u_E^*)^n=0$ on $E^0$. And by Proposition \ref{p13.5}, $(\triangle u_E^*)^n=0$ on $\Omega\backslash \overline{E}$. So $(\triangle u_E^*)^n$ is supported by $\partial E$.
\e{proof}

\b{pro}\label{p13.15}Let $\Omega\subseteq\mathbb{H}^n$ be a strongly pseudoconvex smooth open set. For arbitrary sets $E\Subset\Omega$, we have $C^*(E,\Omega)=\int_\Omega(\triangle u_E^*)^n.$
\e{pro}
\proof First we show that for a compact set $K\subseteq\Omega$,
\beq \label{13.151}C(K,\Omega)=\int_\Omega(\triangle u_K^*)^n=\int_K(\triangle u_K^*)^n.\eeq
The last equality in (\ref{13.151}) follows from Proposition \ref{p13.5} directly. Since $-1\leq u_K^*\leq0$ on $\Omega$, $C(K,\Omega)\geq\int_K(\triangle u_K^*)^n$ by definition. Let $\psi<0$ be a smooth strictly PSH exhaustion function of $\Omega$, we have $A\psi\leq-1$ on $K$ for $A$ large enough.

As in the proof of Proposition \ref{p13.5}, there exists an increasing sequence $\{v_j\}\subseteq PSH(\Omega)$ such that $-1\leq v_j\leq 0$ on $\Omega$, $v_j\leq-1$ on $K$ and $v^*=u_K^*$. We can assume that $v_j\geq A\psi$ on $\Omega$ (otherwise we can replace $v_j$ by $\max\{v_j,A\psi\}$).
Take $\varepsilon\in (0,1)$ and $\omega\in PSH(\Omega)$, $0\leq \omega\leq 1-\varepsilon$. Now we have
$$K\subseteq\{v_j\leq\omega-1\}\subseteq\{A\psi\leq\omega-1\}\subseteq\{A\psi\leq-\varepsilon\}.$$ Note that $v_j\geq A\psi>-\varepsilon >v-1$ near $\partial\Omega$ sufficiently, by the comparison theorem (Lemma \ref{t2.2}) we have
$$\int_K(\triangle \omega)^n\leq \int_{\{v_j\leq\omega-1\}}(\triangle \omega)^n\leq \int_{\{v_j\leq\omega-1\}}(\triangle v_j)^n\leq \int_{\{A\psi\leq-\varepsilon\}}(\triangle v_j)^n.$$ By Lemma \ref{l2}, we have $(\triangle v_j)^n$ converges weakly to $(\triangle u_K^*)^n$  as $j\rightarrow\infty$. Thus
$$\int_K(\triangle \omega)^n\leq\int_{\{A\psi\leq-\varepsilon\}}(\triangle u_K^*)^n=\int_K(\triangle u_K^*)^n.$$
The identity above follows from Proposition \ref{p13.5}. Note that  $$C(K,\Omega)=(1-\varepsilon)^{-n}\sup\{\int_K(\triangle \omega)^n:~\omega\in PSH(\Omega),0\leq \omega\leq1-\varepsilon\}.$$ Then (\ref{13.151}) follows.

Now for every open set $G\Subset\Omega$, we are going to show that
\beq \label{13.152}C^*(G,\Omega)(=C(G,\Omega))=\int_{\overline{G}}(\triangle u_G^*)^n=\int_\Omega(\triangle u_G^*)^n=\int_\Omega(\triangle u_G)^n.
\eeq
Let $K_1\subseteq K_2\subseteq\ldots$ be compact subsets of $G$ with $K_j\subseteq K_{j+1}^0$ and $\bigcup_j K_j=G$. Then $u_{K_j}^*=-1$ on $K_j^0\supseteq K_{j-1}$ and $\lim_ju_{K_j}^*=-1$ on $G$. Since $K_j\subseteq G$, it follows directly from Corollary \ref{c13.5} (1) that $u_G^*\leq u_{K_j}^*$. Then $u_G^*\leq \lim_ju_{K_j}^*$. On the other hand, $\lim_ju_{K_j}^*\leq u_G$ by the definition (\ref{u*}). Therefore $u_G^*\leq \lim_ju_{K_j}^*\leq u_G\leq u_G^*$. Then (\ref{13.152}) follows from (\ref{13.151}) and Proposition \ref{p1.1}(5).

Now let $E\Subset \Omega$ be arbitrary and let $\psi<0$ be a strictly PSH exhaustion function of $\Omega$. For every open set $E\subseteq G\Subset \Omega$, we have $u_G^*\geq A\psi$ and $u_E^*\geq u_G^*$ by Corollary \ref{c13.5} (1) . Note that $0\geq u_G^*\geq A\psi$, $u_G^*(q)\rightarrow0$ as $q\rightarrow\partial\Omega$. By the comparison theorem (Lemma \ref{t2.2}) we have
$$\int_\Omega(\triangle u_E^*)^n\leq \int_\Omega(\triangle u_G^*)^n=C(G,\Omega).$$ Thus $\int_\Omega(\triangle u_E^*)^n\leq C^*(E,\Omega)$ by the definition (\ref{cap*}).

On the other hand, Lemma \ref{choquet} shows that there exists a sequence $\{v_j\}\subseteq PSH$ with $-1\leq v_j\leq 0$, $v_j\geq A\psi$ and $\lim v_j=u_E$ a.e. in $\Omega$. Consider open sets $G_j=\{q\in\Omega,(1+\frac{1}{j})v_j(q)<-1\}$. We have $G_j\supseteq E$, $G_j$ is decreasing and $(1+\frac{1}{j})v_j\leq u_{G_j}\leq u_{G_j}^*$. Note that $u_{G_j}^*\nearrow u_E^*$, we have $C^*(E,\Omega)\leq \lim_j C(G_j,\Omega)=\lim_j \int_\Omega(\triangle u_{G_j}^*)^n=\int_\Omega(\triangle u_{E}^*)^n$. The proposition is proved finally.
\endproof 

\b{cor}\label{c3.3} If $\Omega\in\mathbb{H}^n$ is open, then for each $\omega\Subset\Omega$ and $u\in PSH(\Omega)$, $$\lim_{j\rightarrow\infty}C(\{u<-j\}\cap\omega,\Omega)=0.$$
\e{cor}
\proof   Assume that $u<0$ on $\omega$. Cover $\overline{\omega}$ by a finite union of balls in $\Omega$. By Proposition \ref{p1.1}, we can assume that $\Omega$ is strongly pseudoconvex. Set $P_j=\{u<-j\}\cap\omega$. Then we have $\max\{\frac{u}{j},-1\}\leq u_{P_j}\leq0$. Therefore $\lim_{j\rightarrow\infty}u_{P_j}=0$ a.e. in $\Omega$. By (\ref{13.152}) and Lemma \ref{l2}, we have $\lim_{j\rightarrow\infty}C(P_j,\Omega)=0.$
\endproof

\b{rem}By Corollary \ref{c3.3} above, we can generalize the quasicontinuity theorem (Lemma \ref{t2.4}) to the unbounded case, i.e. for arbitrary $u\in PSH(\Omega) $, for any $\varepsilon>0$, there exists an open subset $\omega\subseteq\Omega$ with $C(\omega)<\varepsilon$ such that $u$ is continuous on $\Omega\backslash\omega$.
\e{rem}

\begin{cor}\label{c3.2} Let $u,v\in PSH\cap L_{loc}^\infty(\Omega)$. If $\limsup_{\zeta\leftarrow \partial\Omega}|u(\zeta)-v(\zeta)|=0,$ and $(\triangle u)^n=(\triangle v)^n$ in $\Omega$, then $u\equiv v$ in $\Omega$.
\end{cor}
\proof It suffices to prove that $u\geq v$. Let $\varphi<0$ be a smooth strictly PSH function in $\Omega$. Suppose that $\{u<v\}$ is not empty, then the set $S=\{u<v+\varepsilon \varphi\}$ is also non-empty for some proper $\varepsilon>0$. Since $u$ and $v+\varepsilon \varphi$ are both subharmonic, by the classical result for the subharmonic functions we know that the set $S$ must have positive Lebesgue measure. By Lemma \ref{t2.2},
$$\int_S(\triangle u)^n\geq \int_S(\triangle (v+\varepsilon\varphi))^n\geq  \int_S(\triangle v)^n+\varepsilon^n \int_S(\triangle \varphi)^n.$$ The last integral over $S$ is strictly positive, so we get a contradiction.
\qed

By Proposition \ref{p13.15} and Corollary \ref{c3.2}, we have
\b{cor}\label{c13.17} Let $\Omega$ be a strongly pseudoconvex smooth open set in $\mathbb{H}^n$ and let $E\Subset\Omega$. Then $C^*(E,\Omega)=0$ if and only if $u_E^*=0$.
\e{cor}
\b{pro}\label{p13.11} Let $\Omega$ be a connected open set in $\mathbb{H}^n$, and let $E\subseteq \Omega$. The following statements are equivalent:\\
(1). $u_E^*\equiv0$;\\
(2). there exists $v\in PSH(\Omega)$, $v\leq0$ such that $E\subseteq \{v=-\infty\}$.
\e{pro}
\proof The implication $(2)\Rightarrow(1)$ is obvious. If $v$ is as in (2), then for each $\varepsilon>0$, $\varepsilon v\leq u_E$, so $u_E=0$ on $\Omega\backslash\{v=-\infty\}$. So $u_E^*\equiv0$.

Now assume that $u_E^*\equiv0$. By Lemma \ref{choquet}, there exists a sequence $\{v_j\}\subseteq PSH(\Omega)$, $-1\leq v_j\leq u_E$, converging increasingly a.e. in $\Omega$ to $u_E^*$. We can extract a subsequence in such a way that $\int_\Omega|v_j|d\lambda<2^{-j}$. Since $v_j\leq0$ and $v_j\leq-1$ on $E$, the function $v:=\sum v_j\leq0,$ and $v=-\infty$ on $E$. As $v$ is the limit of the decreasing sequence of its partial sums and $v\not\equiv-\infty$ in $\Omega$, we have $v\in PSH(\Omega)$ by Proposition \ref{p2.1}. \endproof

Now we prove the Josefson's theorem on $\mathbb{H}^n$, following the proof given in \cite{bed} in pluripotential theory on $\mathbb{C}^n$.

\emph{Proof of Theorem \ref{josefson}.}
 By the definition of locally Q-polar, we can find sets $P_j,\Omega_j$ with $\Omega_j$ strongly pseudoconvex smooth open, $P_j\Subset\Omega_j\Subset\mathbb{H}^n$, $\bigcup_{j\geq1}P_j=P$ and $P_j$ is contained in the $-\infty$ poles of a single plurisubharmonic function in $\Omega_j$. By Proposition \ref{p13.11} and Proposition \ref{p13.15}, we have $C^*(P_j,\Omega_j)=0$.

 Let $i_1,i_2,\ldots$ be a listing of the positive integers such that each one appears infinitely many times. For a sequence $c_1<c_2<\ldots$, $c_j\rightarrow+\infty$, set $B_j=\{q\in \mathbb{H}^n, \|q\|<c_j\}$. We can choose $c_j$ large enough that $\Omega_{i_j}\Subset B_j$ and $|q|-c_j<-1$ on $P_{ij}$. It follows from $C^*(P_{i_j},\Omega_{i_j})=0$ and $\Omega_{i_j}\Subset B_j$ that $C^*(P_{i_j},B_j)=0$. Hence by Corollary \ref{c13.17} the extremal function $u_{P_{ij}}^*$ in $B_j$ is zero and we can find $v_j\in PSH(B_j)$ with $v_j\leq0$ on $B_j$, $v_j\leq-1$ on $P_{i_j}$ and $\int_{B_j}|v_j|dV<2^{-j}$. Define $$\widetilde{v_j}(q)=\left\{
                                                                                                                                      \begin{array}{ll}
                                                                                                                                        \|q\|-c_j, &q\in \mathbb{H}^n\backslash B_j, \\
                                                                                                                                        \max\{v_j(q), \|q\|-c_j\}, & q\in B_j.
                                                                                                                                      \end{array}
                                                                                                                                    \right.$$
Then $\widetilde{v_j}\leq-1$ on $P_{i_j}$ and $\widetilde{v_j}\in PSH(\mathbb{H}^n)$ by Proposition \ref{p2.1}. As $\widetilde{v_j}<0$ on $B_j$ and $\int_{B_j}|v_j|dV<2^{-j}$, $v=\sum_{j=1}^\infty\widetilde{v_j}$ is a PSH function on $\mathbb{H}^n$. Since $\widetilde{v_j}=-1$ on $P_{i_j}$ and each $P_i$ repeated infinitely many times, therefore $v=-\infty$ on $P=\bigcup_{j\geq1}P_j$. This completes the proof.\qed

\b{lem}\label{p13.12} Let $\Omega$ be an open subset of $\mathbb{H}^n$ and let $K_1\supseteq K_2\supseteq\ldots$, $K=\bigcap_j K_j$ be compact subsets of $\Omega$. Then\\
(a) $\lim C(K_j,\Omega)=C(K,\Omega).$\\
(b) $C^*(K,\Omega)=C(K,\Omega)$.\\
In particular, $C^*(K,\Omega)=C(K,\Omega)$ for any compact set $K\subseteq\Omega$.
\e{lem}
\proof (a) follows from Lemma \ref{l2} and Proposition \ref{p13.15}. Note that $K_j$'s are neighborhoods of $K$, (b) follows directly from (a).\endproof

As an application of quasicontinuity theorem, we can prove an interesting inequality for the quaternionic Monge-Amp\`ere operator. Here we follow the proof of the complex case in Demailly \cite{Demailly1991}.
\b{pro}\label{p2}(Demailly's inequality)Let $u,v$ be locally bounded PSH functions on $\Omega$. Then we have an inequality of quaternionic Monge-Amp\`ere measures
\beq (\triangle\max\{u,v\})^n\geq \chi_{\{u\geq v\}} (\triangle u)^n+\chi_{\{u< v\}} (\triangle v)^n.
\eeq
\e{pro}
\b{proof} By changing the roles of $u$ and $v$, it suffices to prove that
\beq\label{11.9}\int_K(\triangle\max\{u,v\})^n\geq\int_K  (\triangle u)^n\eeq
 for every compact set $K\subseteq\{u\geq v\}$. Since $u,v$ are bounded, we may assume that $0\leq u,v\leq 1$ and $0\leq u_\varepsilon,v_\varepsilon\leq 1$, where $u_\varepsilon:=u*\rho_\varepsilon$ is the regularization of $u$. By Lemma \ref{t2.4}, we can assume that $G\subseteq \Omega$ is an open set of small capacity such that $u,v$ are continuous on $\Omega\backslash G$. Then $u_\varepsilon,v_\varepsilon$ converge uniformly to $u,v$ on $\Omega\backslash G$ respectively as $\varepsilon$ tends to $0$. For any $\delta>0$, we can find an arbitrarily small neighborhood $L$ of $K$ such that $u_\varepsilon>v_\varepsilon-\delta$ on $L\backslash G$ for $\varepsilon$ sufficiently small. By Lemma \ref{l2.3}, $(\triangle u_\varepsilon)^n$ converges weakly to $(\triangle u)^n$, so we have
\begin{equation*}\begin{aligned}\int_K  (\triangle u)^n\leq \liminf_{\varepsilon\rightarrow0}\int_L  (\triangle u_\varepsilon)^n\leq&  \liminf_{\varepsilon\rightarrow0}\left(\int_G  (\triangle u_\varepsilon)^n+\int_{L\backslash G}  (\triangle u_\varepsilon)^n\right)\\ \leq& C(G,\Omega)+\liminf_{\varepsilon\rightarrow0}\int_{L\backslash G}  (\triangle \max\{u_\varepsilon+\delta,v_\varepsilon\})^n\\=&C(G,\Omega)+\int_{L\backslash G}  (\triangle \max\{u+\delta,v\})^n
\end{aligned}\end{equation*}
The third inequality above follows from the definition of capacity and the fact that $\max\{u_\varepsilon+\delta,v_\varepsilon\}=u_\varepsilon+\delta$ on a neighborhood of $L\backslash G$.

By taking $L$ very close to $K$ and $C(G,\Omega)$ arbitrarily small, we have $\int_K  (\triangle u)^n\leq \int_{K}  (\triangle \max\{u+\delta,v\})^n.$ Let $\delta\rightarrow0$ to get ($\ref{11.9}$).
\e{proof}

\b{pro}\label{p13.6} Every negligible set $N\subseteq \Omega$ satisfies $C^*(N,\Omega)=0$.
\e{pro}
\proof By Lemma \ref{choquet}, every neglibible set is contained in a Borel negligible set $N=\{v<v^*\}$ with $v=\sup v_\alpha$, and $\{v_\alpha\}$ is an increasing sequence of PSH functions with $v_\alpha\geq-1$ for all $\alpha$. By Lemma \ref{t2.4}, there exists an open set $G\subseteq\Omega$ such that all $v_\alpha$ and $v^*$ are continuous on $\Omega\backslash G$ and $C(G,\Omega)<\varepsilon$. As $G$ is open, $C^*(G,\Omega)=C(G,\Omega)<\varepsilon$.

Write $$N\subseteq G\cup(N\cap(\Omega\backslash G))=G\cup(\bigcup_{\delta,\lambda,\mu}K_{\delta\lambda\mu})$$ with $$K_{\delta\lambda\mu}=\{q\in \overline{\Omega_\delta}\backslash G, v(q)\leq \lambda<\mu\leq v^*(q)\},~ \lambda<\mu,~\lambda,\mu\in \mathbb{Q},~\delta>0.$$
Set $K=K_{\delta\lambda\mu}$ for short. As $v^*$ is continuous and $v$ lower semi-continuous on $\Omega\backslash G$, we see that $K$ is either compact or empty. It suffices to show that $C(K,\Omega)=0$. Take an open set $\omega\Subset \Omega$, we may assume that $v^*\leq0$ on $\overline{\omega}$. Set $\lambda=-1$. Then $v_\alpha\leq0$ on $\omega$ and $v_\alpha\leq v\leq-1$ on $K$. So $v\leq u_K, v^*\leq u_K^*$. And $u_K^*\geq\mu>-1$ on $K$. By Proposition \ref{p2} we have
$$C(K,\omega)=\int_K(\triangle u_K^*)^n\leq\int_K(\triangle\max\{ u_K^*,\mu\})^n\leq|\mu|^nC(K,\omega)$$ as $-1\leq|\mu|^{-1}\max\{ u_K^*,\mu\}\leq0$. Since $|\mu|<1$, we have $C(K,\omega)=0$. By Lemma \ref{p13.12} $C^*(K,\Omega)=C(K,\Omega)=0$. So $C^*(N,\Omega)<\varepsilon$ for every $\varepsilon>0$.
\endproof

\b{cor}\label{c14.4} Let $\Omega\subseteq\mathbb{H}^n$ and $P\subseteq\Omega$. Then $P$ is Q-polar in $\Omega$ if and only if $C^*(P,\Omega)=0$.
\e{cor}
\proof If $C^*(P,\Omega)=0$, then $C^*(P\cap\omega',\omega)=0$ for all concentric balls $\omega'\Subset\omega\Subset\Omega$. By Corollary \ref{c13.17}, $u_{P\cap\omega'}^*=0$. And by Proposition \ref{p13.11}, there exists $0\geq v\in PSH(P\cap\omega')$ such that $P\cap\omega'\subseteq\{v=-\infty\}$, i.e. $P$ is locally Q-polar. Then the conclusion follows from Theorem \ref{josefson} and Corollary \ref{c2.1}.
\endproof

Finally, we get ready to combine all results obtained together to complete the proof of Theorem \ref{t1}.

\vspace{3mm}

 \emph{Proof of Theorem \ref{t1}.} By Corollary \ref{c14.4}, it suffices to show that negligible sets are the same as Q-polar sets. Proposition \ref{p13.6} and Corollary \ref{c14.4} imply that every negligible set is Q-polar. Proposition \ref{p13.2} applied in $\mathbb{H}^n$ shows that each globally Q-polar set is negligible. Then Theorem \ref{t1} follows from Theorem \ref{josefson} and Corollary \ref{c13.17}.
\qed

\b{thm}\label{t3.1}If $\Omega\subseteq\mathbb{H}^n$ is strongly pseudoconvex smooth open, then the function $E\rightarrow C^*(E,\Omega)$ is a generalized capacity. This means:\\ (1). $C^*(\varnothing,\Omega)=0$.\\
(2). If $K_1\supseteq K_2\supseteq\ldots$ is a sequence of compact subsets of $\Omega$, then $$\lim_{j\rightarrow\infty}C^*(K_j,\Omega)=C^*(\bigcap_{j=1}^\infty K_j,\Omega).$$
(3). If $E_1\subseteq E_2\subseteq\ldots$ is a sequence of arbitrary subsets of $\Omega$, then $$\lim_{j\rightarrow\infty}C^*(E_j,\Omega)=C^*(\bigcup_{j=1}^\infty E_j,\Omega).$$
All Suslin (in particular all Borel) subsets $E$ of $\Omega$ are capacitable, that is, $C^*(E,\Omega)=C(E,\Omega)$.
\e{thm}
\proof The property (1) is obvious and (2) was shown in Lemma \ref{p13.12}. By Proposition \ref{p1.1}, for each $E_j$, $C^*(E_j,\Omega)\leq C^*(\bigcup_{j=1}^\infty E_j,\Omega)$. It follows that $\lim_{j\rightarrow\infty}C^*(E_j,\Omega)\leq C^*(\bigcup_{j=1}^\infty E_j,\Omega).$

To prove the opposite inequality, it suffices to show this under the hypothesis that the sets $E_j\Subset\Omega$. Take $\varepsilon,\delta\in(0,1)$. By Theorem \ref{t1}, the sets $\widetilde{E}_j:=\{q\in E_j,u_{E_j}^*>-1\}$ are Q-polar sets and satisfy $C^*(\widetilde{E}_j,\Omega)=0$. By Proposition \ref{p1.1}, the union $F:=\bigcup_{j=1}^\infty\widetilde{E}_j$ satisfies $C^*(F,\Omega)=0$. It follows from the definition (\ref{cap*}) that there exists an open set $G$, $F\subseteq G\subseteq\Omega$, such that $C^*(G,\Omega)<\varepsilon$. Define $$U_j=\{q\in\Omega:u_{E_j}^*<-1+\delta\} ~~\text{and}~~ V_j=U_j\cup G.$$ Since $u_{E_j}^*$ is upper semicontinuous, $U_j,V_j$ are open. Note that $\frac{1}{1-\delta}u_{E_j}^*\leq u_{U_j}^*$ in $\Omega$. Then by the subadditivity of capacity $C(\cdot,\Omega)$ and Proposition \ref{p13.15}, we have\begin{equation*}\begin{aligned}C^*(V_j,\Omega)&\leq C^*(G,\Omega)+C^*(U_j,\Omega)\leq \varepsilon+\int_\Omega(\triangle u_{U_j}^*)^n\\
&\leq \varepsilon+(1-\delta)^{-n}\int_\Omega(\triangle u_{E_j}^*)^n=\varepsilon+(1-\delta)^{-n}C^*(E_j,\Omega).\end{aligned}\end{equation*}
Therefore we have \begin{equation*}\begin{aligned}C^*(\bigcup_{j=1}^\infty E_j,\Omega)&\leq C^*(\bigcup_{j=1}^\infty V_j,\Omega)=C(\bigcup_{j=1}^\infty V_j,\Omega)\\&=\lim_{j\rightarrow\infty}C(V_j,\Omega)\leq \varepsilon+(1-\delta)^{-n}\lim_{j\rightarrow\infty}C^*(E_j,\Omega).\end{aligned}\end{equation*}
Finally let $\varepsilon,\delta\rightarrow0$ to obtain the required estimate. It is the classical Choquet's result that the last conclusion of the theorem follows from (1)-(3). See for example Chapter 2 in \cite{function}, Chapter 3 in \cite{homander}.
\qed

\begin{Acknw}
This work is supported by Natural Science Foundation of SZU (grant no. 201424) and National Nature Science Foundation in China (No. 11401390; No.
11171298).
\end{Acknw}

  \bibliographystyle{elsarticle-num}
 \bibliography{mybibfile}

\end{document}